\documentclass{article}
\usepackage{amssymb}
\usepackage{epsfig}
\usepackage{amsmath}
\usepackage{amsfonts}
\usepackage{tikz}
\usepackage{float}
\usepackage{color}
\usepackage{hyperref}

\newcommand{\R}{{{\Bbb R}}}

\newcommand{\N}{{{\Bbb N}}}

\newtheorem{theorem}{\sc Theorem}[section]
\newtheorem{proposition}{\sc Proposition}[section]
\newtheorem{lemma}{\sc Lemma}[section]
\newtheorem{definition}{\sc Definition}[section]
\newtheorem{remark}{\sc Remark}[section]
\newtheorem{corollary}{\sc Corollary}[section]
\newtheorem{example}{\sc Example}[section]

\title{The method of coupled fixed points and coupled quasisolutions when working with ODE's with arguments of bounded variation\footnote{Partially
supported by FEDER and
Ministerio de Educaci\'on y Ciencia, Spain,
project MTM2010-15314.}}

\author{Rub\'en Figueroa}

\date{}

\begin{document}
\maketitle

\begin{center}
 Departamento de
An\'alise Matem\'atica\\
Facultade de Matem\'aticas\\Universidade de Santiago de Compostela, Campus Vida \\
15782 Santiago de
Compostela, Spain\\
  {\bf e-mail:} ruben.figueroa.sestelo@gmail.com
\end{center}

\def\tanh#1{\,{\normalsize tanh}{\,#1}\,}
\def\simbolo#1#2{#1\dotfill{#2}}
\def\qed{\hbox to 0pt{}\hfill$\rlap{$\sqcap$}\sqcup$\medbreak}
\def\theequation{\arabic{section}.\arabic{equation}}
\def\thesection {\arabic{section}}

\begin{abstract} The aim of this paper is to show the use of the coupled quasisolutions method as a useful technique when treating with ordinary differential equations with functional arguments of bounded variation. We will do this by looking for solutions for a first--order ordinary differential equation with an advanced argument of bounded variation. The main trick is to use the Jordan decomposition of this argument in a nondecreasing part and a nonincreasing one. As a necessary step, we will also talk about coupled fixed points of multivalued operators.
\end{abstract}


\section{Introduction}

In the paper \cite{fp} we proved a new result on the existence of coupled fixed points for mul\-ti\-va\-lued operators, and then we used it to guarantee the existence of coupled quasisolutions and solutions to a certain first--order ordinary differential equation with state--dependent delay. In that paper, the nonlinearity was allowed to have both nondecreasing and nonincreasing arguments and the existence of solutions was obtained under strong Lipschitz conditions. We pointed out there that this tool could be useful when working with arguments of bounded variation, but no literature about this was written  since then. So, the main goal in the present paper is to develop the application of the coupled quasisolutions technique in the framework of arguments of bounded variation, and we do it in an appropiate way, in order to take advantage of the Jordan decompostion and avoid the use of strong assumptions, as Lipschitz--continuity. \\

To show the application of this technique, we will study throughout this paper the existence of solutions for the following first--order problem:

\begin{equation}\label{p}
\left\{
\begin{array}{ll}
x'(t)=f(t,x(t),x(\tau(t))) \ \mbox{ for almost all (a.a.) $t \in I=[a,b]$,} \\
\\
x(t)=\phi(t) \ \mbox{ for all $t \in [b,b+r]$},
\end{array}
\right.
\end{equation}
where, $r \ge 0$, $\tau$ is a measurable function such that $\tau(t) \ge t$ for a.a. $t$, that is, $\tau$ is an advanced argument, and $\phi$ is a bounded function which represents the final state of the solution. By a solution of (\ref{p}) we mean a function $x \in \mathcal{C}[a,b+r]$ such that $x_{|I} \in AC(I)$ and $x$ satisfies both the differential equation (almost everywhere on $I$) and the final condition. We refer the readers to papers \cite{dykjan}, \cite{fig}, \cite{jankowski} to see more results on the existence of solutions and some applications of first--order problems with advance. \\

This paper is organized as follows: in Section $2$ we gather some preliminary concepts and results involving functions of bounded variation and coupled fixed points of mul\-ti\-va\-lued operators. These preliminaries are used later, in Section $3$, to prove the existence of quasisolutions and solutions for problem (\ref{p}). Finally, in Section $4$ an example of application is provided too.

\section{Preliminaries on bounded variation and coupled fixed points of multivalued operators}

In this section we introduce some preliminaries that we will uses throughout this work. First, we remember some concepts about functions of bounded variation. The reader can see more about this in the monographs \cite{apostol}, \cite{hest}.

\begin{definition} Given a function $f:I=[a,b] \subset \R \longrightarrow \R$, and a partition $P=\{x_0, \ldots, x_n\}$ of $I$, we define the variation of $f$ relative to the partition $P$ as the number
$$
V(f,P)=\sum_{i=1}^n |f(x_{i})-f(x_{i-1})|,$$
and we define the total variation of $f$ on $I$ as
$$
V_a^b(f)=\sup_{P \in \mathcal{P}} V(f,P),$$
where $\mathcal{P}=\{P \, : \, P \mbox{ is a partition of $I$}\}$. \\
We say that $f$ is a function of bounded variation on $I$ if $V_a^b(f) <+\infty$. In that case, we write $f \in BV(I)$.
\end{definition}

Functions of bounded variation satisfy the following well--known result, which becomes essential now for our porpuses.

\begin{proposition}\label{jordan}{\bf (Jordan decomposition)} A function $f$ is of bounded variation on $I$ if and only if there exist a nondecreasing function, $g$, and a nonincreasing one, $h$, such that
$$
f(t)=g(t)+h(t) \ \mbox{ for all $t \in I$.}$$
\end{proposition}

The proof of Proposition (\ref{jordan}) uses the fact that the function $t \in I \longrightarrow V_a^t(f)$ is nondecreasing and $t \longrightarrow f(t)-V_a^t(f)$ is nonincreasing, and so the desired decomposition is
$$f(t)=V_a^{t}(f) +f(t)-V_a^{t}(f).$$ We remark that this decomposition is not unique. Finally, notice that, as a consequence of this result, every function of bounded variation is almost--everywhere di\-ffe\-ren\-tia\-ble. \\

The set $BV(I)$ is an algebra which is included neither in the set of continuous functions nor in its complementary. Indeed, if $f$ is monotone on $[a,b]$ then $V_a^b(f)=|f(b)-f(a)|$, and so $f \in BV(I)$. Then, there exists discontinuous functions which are of bounded variation (for example, an step function). In fact, it is also a well--known fact that if $f \in BV(I)$ then $f$ has only ``jump'' discontinuities. On the other hand, there exists continuous functions which are not of bounded variation, as (see \cite[Example 6.3.1]{apostol})
$$
f(t)=\left\{
\begin{array}{ll}
t \cos \left(\dfrac{\pi}{2t}\right) \quad &\mbox{if $0<t \le 1$}, \\
& \\
0 \quad &\mbox{if $t=0$}.
\end{array}
\right.
$$

To obtain our main result, we will use a generalized monotone method in presence of lower and upper solutions. This is a very well--known tool which is extensively used in the literature of ordinary differential equations. The classical version of this technique uses a pair of monotone sequences which will converge to the extremal solutions of the problem. The generalized version of this technique was developed in \cite{hela}, and it is used when the nonlinearity has discontinuous arguments and so the pair of monotone sequences is replaced by a monotone operator. As a novelty which respect to the method developed in \cite{hela} and related references, we will use here a multivalued operator (that is, a set--valued mapping) defined in a product space, and so we will look for coupled fixed points. We concrete this idea in the following lines.

\begin{definition} A metric space $X$ equipped with a partial ordering $\le$ is an ordered metric space if the intervals $[x)=\{y \in X \, :\, x \le y\}$ and $(x]=\{y \in X \, :\, y \le x\}$ are closed for every $x \in X$. Let $P$ be a subset of an ordered metric space; an operator $A:P \times P \longrightarrow P$ is said to be mixed monotone if $A(\cdot,x)$ is nondecreasing and $A(x,\cdot)$ is nonincreasing for each $x \in P$. We say that $A$ satisfies the mixed monotone convergence property (m.m.c.p.) if $(A(v_j,w_j))_{j=1}^{\infty}$ converges in $Y$ whenever $(v_j)_{j=1}^{\infty}$ and $(w_j)_{j=1}^{\infty}$ are sequences in $P$, one being nondecreasing and the other nonincreasing.
\end{definition}

\begin{definition} Let $\overline{X}$ be a subset of an ordered metric space $X$. We define a multivalued operator in the product $\overline{X} \times \overline{X}$ as a mapping
$$
\mathcal{A}: \overline{X} \times \overline{X} \longrightarrow 2^{\overline{X}} \backslash \emptyset.$$
We say that $v,w \in \overline{X}$ are coupled fixed points of $\mathcal{A}$ if $v \in \mathcal{A}(v,w)$ and $w \in \mathcal{A}(w,v)$. We say that $v_*,w^* \in \overline{X}$ are the extremal coupled fixed points of $\mathcal{A}$ in $\overline{X}$ if $v_*,w^*$ are coupled fixed points of $\mathcal{A}$ and if $v,w \in \overline{X}$ is another pair of coupled fixed points of $\mathcal{A}$ then $v_* \le v$ and $w \le w^*$.
\end{definition}

\begin{theorem}{\bf \cite[Theorem 2.1]{fp}}\label{multi} Let $Y$ be a subset of an ordered me\-tric space $X$, $[\alpha,\beta]$ a nonempty closed interval in $Y$ and $\mathcal{A}:[\alpha,\beta] \times [\alpha,\beta] \longrightarrow 2^{[\alpha,\beta]} \backslash \emptyset$ a multivalued operator.

If for all $v,w \in [\alpha,\beta]$ there exist
$$A_*(v,w)=\min \mathcal{A}(v,w) \in [\alpha,\beta], \quad A^*(v,w)=\max \mathcal{A}(v,w) \in [\alpha,\beta],$$ and the (single--valued) operators $A_*$ and $A^*$ are mixed monotone and satisfy the m.m.c.p., then $\mathcal{A}$ has the extremal coupled fixed points in $[\alpha,\beta]$, $v_*,v^*$. Moreover, they satisfy the following characterization:
\begin{equation}\label{minmax}
(v_*,w^*)=\min_{\preceq} \{(v,w) \, : \, (A_*(v,w),A^*(w,v)) \preceq (v,w)\},
\end{equation}
where $$(v,w) \preceq (\overline{v},\overline{w}) \Longleftrightarrow v\le \overline{v}, w \ge \overline{w}.$$

\end{theorem}

\section{Main result}

Now we develop our generalized monotone method applied to problem (\ref{p}). To do this, throughout this section we will assume the following: \\

$(H_1)$ There exists a closed interval $J \subset \R$ such that for a.a. $t \in I$ and all $x \in \R$ the function $f(t,x,\cdot)$ is of bounded variation on $J$.\\

Assumption $(H_1)$ implies that there exists a nondecreasing function, $g$, and a nonincreasing one, $h$, such that
$$
f(t,x,\cdot)=g(t,x,\cdot)+h(t,x,\cdot)$$
for all $(t,x) \in I \times \R$. \\

Now we define what we mean by lower and upper solutions for problem (\ref{p}).

\begin{definition}\label{lu} We say that $\alpha,\beta \in \mathcal{C}[a,b+r]$ are, respectively, a lower and an upper solution for problem (\ref{p}) if $\alpha_{|I},\beta_{|I} \in AC(I)$, $$\left[\min_{t \in [a,b+r]} \alpha(t),\max_{t \in [a,b+r]} \beta(t)\right] \subset J,$$ the compositions
$$
t \longmapsto f(t,\alpha(t),y), \quad t \longmapsto f(t,\beta(t),y)$$
are measurable for all $y \in J$ and the following inequalities hold:
$$
\left\{
\begin{array}{ll}
\alpha'(t) \ge g(t,\alpha(t),\beta(\tau(t))) + h(t,\alpha(t),\alpha(\tau(t))) \ \mbox{ for a.a. $t \in I$,} \\
\\
\alpha(t) \le \phi(t) \ \mbox{ for all $t \in [b,b+r],$}
\end{array}
\right.
$$
$$
\left\{
\begin{array}{ll}
\beta'(t) \le g(t,\beta(t),\alpha(\tau(t))) + h(t,\beta(t),\beta(\tau(t))) \ \mbox{ for a.a. $t \in I$,} \\
\\
\beta(t) \ge \phi(t) \ \mbox{ for all $t \in [b,b+r],$}
\end{array}
\right.
$$
\end{definition}

\begin{remark} Notice that, under the previous definition, the lower and the upper solution appear ``coupled''. On the other hand, it is assumed that
$$
\min_{t \in [a,b+r]} \alpha(t) \le \max_{t \in [a,b+r]} \beta(t).$$
This is not an strong assumption, taking into account that, as usual, we will ask the lower and the upper solution to be well--ordered in the whole interval $[a,b+r]$. \\

On the other hand, the fact of being $t \longmapsto f(t,\alpha(t),y)$ and $t \in I \longmapsto f(t,\beta(t),y)$ measurable for all $y \in J$ implies that the compositions
$$
t \in I \longmapsto g(t,\alpha(t),\beta(\tau(t))) + h(t,\alpha(t),\alpha(\tau(t))),$$
and
$$
t \in I \longmapsto g(t,\beta(t),\alpha(\tau(t))) + h(t,\beta(t),\beta(\tau(t)))$$
are measurable too, because of being $g$ and $h$ monotone with respect to their last variables.
\end{remark}

As we said in Introduction, an essential tool in our work refers the use of coupled quasisolutions. So, we introduce now this concept.

\begin{definition}\label{qs} We say that two functions $x_*,x^* \in \mathcal{C}[a,b+r]$ are coupled quasisolutions of problem (\ref{p}) if $x_{*|I},x^*_{|I} \in AC(I)$, $x_*(t)=x^*(t)=\phi(t)$ for all $t \in [b,b+r]$ and for a.a. $t \in I$ they satisfy
$$
\left\{
\begin{array}{ll}
x_*(t)=g(t,x_*(t),x^*(\tau(t)))+h(t,x_*(t),x_*(\tau(t))), \\
\\
x^*(t)=g(t,x^*(t),x_*(\tau(t)))+h(t,x^*(t),x^*(\tau(t))).
\end{array}
\right.
$$
We say that these coupled quasisolutions are extremal in a subset $\overline{X} \subset \mathcal{C}[a,b+r]$ if $x_*,x^* \in \overline{X}$ and $x_*(t) \le x_1(t)$, $x_2(t) \le x^*(t)$ whenever that $x_1,x_2 \in \overline{X}$ is another pair of quasisolutions.
\end{definition}

As an auxiliar tool for proving our main result, we need the following maximum principle related to problems with advance. Compare it with \cite[Lemma 3.2]{fig}, \cite[Lemma 1]{jankowski}.

\begin{lemma}\label{mp} Let $\tau:I \longrightarrow [a,b+r]$ be a measurable function such that $\tau(t) \ge t$ for a.a. $t \in I$ and assume that $p \in \mathcal{C}[a,b+r]$ is such that $p_{|I}  \in AC(I)$ and satisfies
$$
\left\{
\begin{array}{ll}
p'(t) \ge K(t) p(t) - L(t) p(\tau(t)) \ \mbox{ for a.a. $t \in I$}, \\
\\
p(t)=0 \ \mbox{ for all $t \in [b,b+r]$},
\end{array}
\right.
$$
where $K,L \in L^1(I)$ and $L \ge 0$ a.e. \\
If
\begin{equation}\label{int}
\int_a^b (K_-(t)+L(t)) \, dt <1,\end{equation}
where $K_-=\max\{-K,0\}$, then $p(t) \le 0$ for all $t \in [a,b+r]$.
\end{lemma}

\noindent {\bf Proof.} Let $t_1 \in [a,b+r]$ such that
$$
p(t_1)=\max_{t \in [a,b+r]} p(t)$$
and assume by contradiction that $p(t_1) >0$. Then, $t_1 \in [a,b)$. Now, let $t_2 \in (t_1,b]$ such that $p(t_2)=0$ and $p(t) \ge 0$ for all $t \in [t_1,t_2]$. Now, integrating between $t_1$ and $t_2$ we obtain
$$
p(t_1)=-\int_{t_1}^{t_2} p'(t) \, dt \le -\int_{t_1}^{t_2} K(t) p(t) \, dt + \int_{t_1}^{t_2} L(t)p(\tau(t))) \, dt \le p(t_1) \int_{t_1}^{t_2} (K_-(t) + L(t)) \, dt,$$
and then condition (\ref{int}) provides the contradiction $p(t_1) < p(t_1)$. \qed

The main result on this paper concerns the existence of extremal quasisolutions and solutions for problem (\ref{p}). It is the following.

\begin{theorem}\label{main} Assume $(H_1)$ and that there exists $\alpha,\beta \in \mathcal{C}[a,b+r]$ which are, respectively, a lower and an upper solution for problem (\ref{p}) such that $\alpha(t) \le \beta(t)$ for all $t \in [a,b+r]$ and
$$
E=\left[\min_{t \in [a,b+r]} \alpha(t),\max_{t \in [a,b+r]} \beta(t)\right] \subset J.$$
Assume moreover that the following conditions hold:

\begin{enumerate}

\item[$(H_2)$] For each $\gamma_1,\gamma_2 \in [\alpha,\beta]=\{\gamma \in \mathcal{C}[a,b+r] \, : \, \alpha(t) \le \gamma(t) \le \beta(t) \ \mbox{ for all $t \in [a,b+r]$}\}$ the final value problem
\begin{equation*}
(P_{\gamma_1,\gamma_2})\left\{
\begin{array}{ll}
x'(t)=F_{\gamma_2,\gamma_1}(t,x(t)):=g(t,x(t),\gamma_2(\tau(t)))+h(t,x(t),\gamma_1(\tau(t))) \ \mbox{ for a.a. $t \in I$}, \\
\\
x(b)=\phi(b)
\end{array}
\right.
\end{equation*}
has the extremal solutions in $[\alpha,\beta]$;

\item[$(H_3)$] There exists $\psi \in L^1(I,[0,+\infty))$ such that for a.a. $t \in I$, all $x \in [\alpha(t),\beta(t)]$ and all $y_1,y_2 \in [\alpha(\tau(t)),\beta(\tau(t))]$ we have
    $$
    |g(t,x,y_1)+h(t,x,y_2)| \le \psi(t);$$

\item[$(H_4)$] There exists $K_1,K_2,L_1,L_2 \in L^1(I)$ such that $L_1,L_2 \ge 0$ a.e. and
$$
g(t,\overline{x},y)-g(t,x,\overline{y}) \ge K_1(t)(\overline{x}-x)-L_1(t)(\overline{y}-y),$$
$$
h(t,\overline{x},\overline{y})-h(t,x,y) \ge K_2(t)(\overline{x}-x)-L_2(t)(\overline{y}-y)$$
whenever that $\alpha(t) \le x \le \overline{x} \le \beta(t)$ and $$\min_{s \in [b,b+r]} \phi(s)-\int_t^b \psi(s) \, ds \le y \le \overline{y} \le \max_{s \in [b,b+r]} \phi(s) + \int_t^b \psi(s) \, ds.$$ \\
Moreover,
\begin{equation}\label{int2}
\int_a^b (K_-(t) + L(t)) \, dt <1,
\end{equation}
where $K=K_1 + K_2$, $L=L_1 + L_2$ and $K_-(t)=\max\{-K(t),0\}$.

\end{enumerate}

In these conditions, problem (\ref{p}) has a unique solution in $[\alpha,\beta]$.

\end{theorem}

\noindent {\bf Proof.} We consider the space $X=\mathcal{C}[a,b+r]$ endowed with the ordering
$$
\gamma_1 \le \gamma_2 \Longleftrightarrow \gamma_1(t) \le \gamma_2(t) \ \mbox{ for all $t \in [a,b+r]$},$$
and we define a multivalued operator
$$
\mathcal{A}:[\alpha,\beta] \times [\alpha,\beta] \subset X \times X \longrightarrow 2^{[\alpha,\beta]} \backslash \emptyset$$
as follows: for each $\gamma_1,\gamma_2 \in [\alpha,\beta]$ we have $x \in \mathcal{A}(\gamma_1,\gamma_2)$ if and only if $x \in [\alpha,\beta]$, $x_{I}$ is a solution of $(P_{\gamma_1,\gamma_2})$ and $x_{|[b,b+r]} = \phi$.\\

\noindent {\it Step $1$: Operator $\mathcal{A}$ has the extremal coupled fixed points in $[\alpha,\beta]$. } By virtue of condition $(H_2)$, operator $\mathcal{A}$ is well--defined and there exist
$$
A_*=\min \mathcal{A}(\gamma_1,\gamma_2), \quad A^*=\max \mathcal{A}(\gamma_1,\gamma_2).$$
We will show now that $A_*,A^*$ are mixed monotone and satisfy m.m.c.p. So, let $$\gamma_1,\overline{\gamma}_1,\gamma_2,\overline{\gamma}_2 \in [\alpha,\beta]$$ such that $\gamma_1 \le \overline{\gamma}_1$, $\gamma_2 \le \overline{\gamma}_2$ and put $$x_1=A_*(\gamma_1,\gamma_2),\quad \overline{x}_1=A_*(\overline{\gamma}_1,\gamma_2),\quad \overline{x}_2=A_*(\gamma_1,\overline{\gamma}_2).$$ Then for all $t \in [b,b+r]$ we have that $x_1(t)=\overline{x}_1(t)=\overline{x}_2(t)$ and for a.a. $t \in I$ we have
$$
\overline{x}_1'(t)=g(t,\overline{x}_1(t),\gamma_2(\tau(t)))+h(t,\overline{x}_1,\overline{\gamma}_1(\tau(t))) \le g(t,\overline{x}_1(t),\gamma_2(\tau(t)))+h(t,\overline{x}_1,\gamma_1(\tau(t))),$$
and so $\overline{x}_1$ is an upper solution for problem $(P_{\gamma_1,\gamma_2})$. The fact of being $x_1$ the least solution of this problem in $[\alpha,\beta]$ implies that $\overline{x}_1 \ge x_1$ and so $A_*(\cdot,\gamma_2)$ is nondecreasing. On the other hand,
$$
x_1'(t)=g(t,x_1(t),\gamma_2(\tau(t))) + h(t,x_1(t),\gamma_1(\tau(t))) \le g(t,x_1(t),\overline{\gamma}_2(\tau(t))) + h(t,x_1(t),\gamma_1(\tau(t))),$$
and so $x_1$ is an upper solution for problem $(P_{\gamma_1,\overline{\gamma}_2})$. Then $x_1 \ge \overline{x}_2$ and so the mapping $A(\gamma_1,\cdot)$ is nonincreasing. In the same way we show that $A^*$ is mixed monotone. \\

To see that $A_*, A^*$ satisfy the m.m.c.p., let $(v_j)_{j=1}^{\infty}$, $(w_j)_{j=1}^{\infty}$ be sequences in $[\alpha,\beta]$, one being nondecreasing and the other nonincreasing. As $A_*,A^*$ are mixed monotone and bounded, we obtain that the sequences $(A_*(v_j,w_j))_{j=1}^{\infty}, (A^*(v_j,w_j))_{j=1}^{\infty}$ have their pointwise limit, say $z_*, z^*$. As $(A_*(v_j,w_j))_{j=1}^{\infty}, (A^*(v_j,w_j))_{j=1}^{\infty}$ are constant in $[b,b+r]$, the convergence is uniform in this interval. On the other hand, for $t,s \in I$, $s <t$, and $j \in \N$ we have
$$
|z_j^*(t)-z_j^*(s)| \le \int_s^t |(g(r,z_j^*(r),w_j(\tau(r))) + h(t,z_j^*(r),v_j(\tau(r))))| \, dr \le \int_t^s \psi(r) \, dr,$$
and so $(z_j^*)_{j=1}^{\infty}$ converges to $z^*$ uniformly on $I$. The same argument is valid for $z_*$. \\

By application of Theorem \ref{multi}, operator $\mathcal{A}$ has the extremal coupled fixed points in $[\alpha,\beta]$, say $x_*,x^*$. \\

\noindent {\it Step $2$: Problem (\ref{p}) has the extremal quasisolutions in $[\alpha,\beta]$.} Indeed, we well show that the extremal coupled fixed points of operator $\mathcal{A}$, $x_*$, $x^*$, correspond with these extremal quasisolutions. First, it is clear that if $x,\overline{x} \in [\alpha,\beta]$ are coupled fixed points of $\mathcal{A}$ then they are coupled quasisolutions of problem (\ref{p}). On the other hand, if $x, \overline{x}$ are quasisolutions of problem (\ref{p}) then $A_*(x,\overline{x}) \le x$ and $A^*(\overline{x},x) \ge \overline{x}$ and so characterization (\ref{minmax}) implies that $x_* \le x$ and $\overline{x} \le x^*$. This shows that $x_*,x^*$ are the extremal quasisolutions of problem (\ref{p}) in $[\alpha,\beta]$. \\

\noindent {\it Step $3$: Problem (\ref{p}) has a unique solution in $[\alpha,\beta]$.} We will prove this by showing that the extremal quasisolutions $x_*,x^*$ are, in fact, the same function, and thus defining a solution of the problem. This solution must be unique in $[\alpha,\beta]$ because if $\overline{x} \in [\alpha,\beta]$ is a solution of (\ref{p}) then the pair $\overline{x},\overline{x}$ is also a quasisolution and then $x_* \le \overline{x} \le x^*$. \\
To see that $x_*=x^*$, first notice that as $(x_*,x^*)$ is a pair of quasisolutions, then the reversed pair, $(x^*,x_*)$, are quasisolutions too, and then extremality implies $x_* \le x^*$. Moreover, condition $(H_3)$ implies that for a.a. $t \in I$
$$
x_*(t),x^*(t) \in \left[ \phi(b) - \int_t^b \psi(s) \, ds,\, \phi(b) + \int_t^b \psi(s) \, ds\right].$$

Now, define the function $p(t)=x^*(t)-x_*(t) \ge 0$. On the one hand, $p(t)=0$ for all $t \in [b,b+r]$. On the other hand, condition $(H_4)$ implies for a.a. $t \in I$ that
\begin{align}
p'(t)&=g(t,x^*(t),x_*(\tau(t)))-g(t,x_*(t),x^*(\tau(t)))+h(t,x^*(t),x^*(\tau(t)))-h(t,x_*(t),x_*(\tau(t))) \\
&\ge K(t)(x^*(t)-x_*(t))-L(t)(x^*(\tau(t))-x_*(\tau(t))),\end{align}
and then by virtue of Lemma \ref{mp} we obtain that $p(t) \le 0$ on $I$. We conclude that $p(t)=0$ for all $t \in [a,b+r]$, that is, $x_*=x^*$. This ends the proof. \qed

\begin{remark} Now we point out some remarks related to Theorem \ref{main}:

\begin{enumerate}
\item{Condition $(H_2)$ could be replaced by any result on the existence of extremal solutions between lower and upper solutions for problem $(P_{\gamma_1,\gamma_2})$. For example, as well--known, if $F_{\gamma_1,\gamma_2}$ is a Carath\'eodory function then $(H_3)$ implies that $(P_{\gamma_1,\gamma_2})$ has the extremal solutions between $\alpha$ and $\beta$. Moreover, there exists a very extensive literature about the existence of extremal solutions for problem $(P_{\gamma_1,\gamma_2})$ for discontinuous $F_{\gamma_1,\gamma_2}$. The reader is referred to \cite{fp}, \cite{hasrzy}, \cite{pouso}, \cite{rlp} and references therein for some results of this type. Notice that although most of these references deal with initial value problems, these results can easily be adapted for final value problems. Finally, notice that $(H_2)$ implies, in particular, measurability of the composition $t \in I \longmapsto F_{\gamma_1,\gamma_2}(t,x(t))$ for all $x \in [\alpha,\beta]$;}

\item{As we said in Section $2$, a function of bounded variation has only ``jump'' dis\-con\-ti\-nui\-ties. Although condition $(H_4)$ implies that for a.a. $t \in I$ the function $f$ is continuous with respect to its third variable in the interval
$$\left[\min_{s \in [b,b+r]} \phi(s)-\int_t^b \psi(s) \, ds, \max_{s \in [b,b+r]} \phi(s) + \int_t^b \psi(s) \, ds\right],$$
a countable number of discontinuities are allowed to exist outside this interval. Moreover, notice that this interval can be improved if we find another function $\widetilde{\psi}$ satisfying $(H_3)$ and such that $\widetilde{\psi}(t) \le \psi(t)$ for almost all $t$;}

\item{For almost all $t \in I$ and all $x \in [\alpha(t),\beta(t)]$ the function $f_{t,x}(\cdot)=f(t,x,\cdot)$ is of bounded variation in $[\alpha(\tau(t)),\beta(\tau(t))]$ and so there exists in this interval a decomposition $f_{t,x}(\cdot)=g_{t,x}(\cdot)+h_{t,x}(\cdot)$, with $g$ nondecreasing and $h$ nonincreasing. Although all conditions in Theorem \ref{main} are stated for an arbitrary Jordan decomposition of this type, all of them can be rewritten with
    $$
    g_{t,x}(y)=V_{A}^{y}(f), \quad h_{t,x}(y)=f_{t,x}(y)-V_{A}^{y}(f),$$
    for any choice of $A \le \min \{\alpha(t) \, : \, t \in [a,b+r]\}$, $A \ge \min J$.}

\end{enumerate}

\end{remark}

Theorem \ref{main} provides, in particular, a new result on the existence of extremal solutions for problem (\ref{p}) in the case that function $f$ is nonincreasing with respect to its third variable. In this case, the nondecreasing part of the Jordan decomposition of $f$ does not exist, and so the lower and upper solutions introduced in Definition \ref{lu} appear uncoupled. Moreover, a pair of quasisolutions in the sense of Definition \ref{qs} become, in fact, a pair of solutions, and so extremal quasisolutions provided by Theorem \ref{main} reduces to extremal solutions. We specify these ideas in the following corollary.

\begin{corollary} Assume that there exist $\alpha,\beta \in \mathcal{C}[a,b+r]$ such that $\alpha_{I},\beta_{I} \in AC(I)$, $\alpha \le \beta$ on $[a,b+r]$ and the following inequalities hold:
$$
\alpha'(t) \ge f(t,\alpha(t),\alpha(\tau(t))) \ \mbox{ for a.a. $t \in I$}, \quad \alpha(t) \le \phi(t) \ \mbox{ for all $t \in [b,b+r]$},$$
$$
\beta'(t) \le f(t,\beta(t),\beta(\tau(t))) \ \mbox{ for a.a. $t \in I$}, \quad \beta(t) \ge \phi(t) \ \mbox{ for all $t \in [b,b+r]$}.$$
Assume moreover that the following conditions hold:

\begin{enumerate}

\item[$(H_2)'$] For all $\gamma \in [\alpha,\beta]$ the final value problem
$$
x'(t)=f(t,x(t),\gamma) \ \mbox{ for a.a. $t \in I$}, \quad x(b)=\phi(b),$$
has the extremal solutions in $[\alpha,\beta]$;

\item[$(H_3)'$] There exists $\psi \in L^1(I,[0,+\infty))$ such that for a.a. $t \in I$, all $x \in [\alpha(t),\beta(t)]$ and all $y \in [\alpha(\tau(t)),\beta(\tau(t))]$ we have
    $$
    |f(t,x,y)| \le \psi(t);$$
\item[$(H_4)'$] For a.a. $t \in I$ and all $x \in [\alpha(t),\beta(t)]$ the function $f(t,x,\cdot)$ is nonincreasing.

\end{enumerate}

In these conditions problem (\ref{p}) has the extremal solutions in $[\alpha,\beta]$.

\end{corollary}

\section{An example of application}

We finish this work with an example of application of our main result.

\begin{example} Consider the following problem with advance:

\begin{equation}\label{ex}
\left\{
\begin{array}{ll}
x'(t)=f(x(4t)) \ \mbox{ for a.a. $t \in I=\left[0,\dfrac{\pi}{8}\right]$},\\
\\
x(t)= \phi(t)=\dfrac{1}{2}\left(x-\dfrac{\pi}{2}\right)\sin\left(\dfrac{1}{x-\pi/2}\right) \ \mbox{ for all $t \in \left[\dfrac{\pi}{8},\dfrac{\pi}{2}\right]$},
\end{array}
\right.
\end{equation}
where $f$ is defined as follows: for each $n\in \{1,2,\ldots\}$ we have
$$
f(y)=\left\{
\begin{array}{ll}
\dfrac{1}{10}y, \quad &\mbox{if $y\in(2n-2,2n-1]$}, \\
& \\
\dfrac{4n-1}{10}-\dfrac{1}{10}y, \quad &\mbox{if $y \in(2n-1,2n]$},
\end{array}
\right.
$$
and for $y \le 0$ we define $f(y)=-f(-y)$. \\

Thus defined, $f$ is a function of bounded variation in any bounded interval of $\R$. Moreover, $f$ has a countable number of both downwards and upwards discontinuities. We will construct later a pair $(\alpha,\beta)$ of coupled lower and upper solutions for problem (\ref{p}) such that for all $t \in \left[0,\dfrac{\pi}{2}\right]$ we have
$$
-\dfrac{\pi}{2} \le \alpha(t) \le \beta(t) \le \dfrac{\pi}{2},$$
and then it suffices to consider a Jordan decomposition of $f$ in the interval $[-2,2]$. So, we put $f=g+h$, with:
$$
g(y)=V_{-2}^{y}(f)=\left\{
\begin{array}{ll}
\dfrac{1}{10}y+\dfrac{2}{10}, \quad &\mbox{if $y \in [-2,-1)$},\\
\\
\dfrac{1}{10}y+\dfrac{3}{10}, \quad &\mbox{if $y \in [-1,1]$}, \\
\\
\dfrac{1}{10}y+\dfrac{4}{10}, \quad &\mbox{if $y \in (1,2]$},
\end{array}
\right.
$$
$$
h(y)=f(y)-V_{-2}^{y}(f)=\left\{
\begin{array}{ll}
-\dfrac{2}{10}y-\dfrac{5}{10}, \quad &\mbox{if $y \in [-2,-1)$},\\
\\
-\dfrac{3}{10}, \quad &\mbox{if $y \in [-1,1]$}, \\
\\
-\dfrac{2}{10}y-\dfrac{1}{10}, \quad &\mbox{if $y \in (1,2]$}.
\end{array}
\right.
$$
\end{example}

We will show now that the functions $\alpha(t)=t-\dfrac{\pi}{2}=-\beta(t)$ are coupled lower and upper solutions for problem (\ref{p}). First, we have that $\alpha(t) \le \phi(t) \le \beta(t)$ for all $t \in \left[\dfrac{\pi}{8},\dfrac{\pi}{2}\right]$. On the other hand, for a.a. $t \in I$ we have the following:

\begin{enumerate}

\item[$(i)$] If $\dfrac{\pi}{2}-4t \in (1,2]$ then $4t - \dfrac{\pi}{2} \in [-2,-1)$ and so
$$
V_{-2}^{\beta(4t)}(f) + \alpha(4t) - V_{-2}^{\alpha(4t)}(f)=\dfrac{8}{10} \left(\dfrac{\pi}{2}-4t\right) = \dfrac{-8}{10} \left(\dfrac{\pi}{2}-4t\right) + \dfrac{2}{10} \le 1=\alpha'(t),$$
$$
V_{-2}^{\alpha(4t)}(f) + \beta(4t) - V_{-2}^{\beta(4t)}(f)= \dfrac{8}{10} \left(\dfrac{\pi}{2}-4t\right) - \dfrac{2}{10} \ge -1=\beta'(t).$$

\item[$(ii)$] If $\dfrac{\pi}{2}-4t \in [0,1]$ then $4t - \dfrac{\pi}{2} \in [-1,0]$ and so
$$
V_{-2}^{\beta(4t)}(f) + \alpha(4t) - V_{-2}^{\alpha(4t)}(f)=\dfrac{8}{10} \left(\dfrac{\pi}{2}-4t\right) = \dfrac{-8}{10} \left(\dfrac{\pi}{2}-4t\right)\le 1=\alpha'(t),$$
$$
V_{-2}^{\alpha(4t)}(f) + \beta(4t) - V_{-2}^{\beta(4t)}(f)= \dfrac{8}{10} \left(\dfrac{\pi}{2}-4t\right) \ge -1=\beta'(t).$$

\end{enumerate}

Then, $\alpha$ and $\beta$ are coupled lower and upper solutions for problem (\ref{ex}), satisfying $\alpha \le \beta$ on $\left[0,\dfrac{\pi}{2}\right]$. \\

Now we check condition $(H_3)$: We have for a.a. $t \in I$, all $x \in [\alpha(t),\beta(t)]$ and all $y_1,y_2 \in [\alpha(4t),\beta(4t)]$ that
$$
|g(t,x,y_1)+h(t,x,y_2)| \le \dfrac{6}{10}+\dfrac{5}{10},$$

and so condition $(H_2)$ is satisfied with $\psi \equiv \dfrac{11}{10}$. \\

Finally, notice that for a.a. $t \in I$ it is
$$
\left[\min_{s \in [\frac{\pi}{8},\frac{\pi}{2}]} \phi(s)-\int_t^{\pi/8} \psi(s) \, ds, \max_{s \in [\frac{\pi}{8},\frac{\pi}{2}]} \phi(s) + \int_t^{\pi/8} \psi(s) \, ds \right] \subset \left[-1,1\right],$$
and so condition $(H_4)$ is satisfied with
$$
L_1 \equiv \dfrac{1}{10}, \quad L_2 \equiv 0.$$

By application of Theorem \ref{main} we conclude that problem (\ref{ex}) has exactly one solution in the functional interval
$$
\left[4t-\dfrac{\pi}{2},\dfrac{\pi}{2}-4t\right].$$

\end{document}